\documentclass{commat}

\usepackage{euscript}
\usepackage[shortlabels]{enumitem}
\usepackage{graphicx}
\usepackage{longtable}

\newcommand\restr[2]{{
		\left.\kern-\nulldelimiterspace 
		#1 
		\right|_{#2} 
}}

\newcommand{\ad}[1]{\mathrm{ad}_{#1}}
\newcommand{\bt}{\mathfrak{b}}
\newcommand{\Co}{\mathcal{C}}
\newcommand{\Cp}{\mathbb{C}}
\newcommand{\epsi}{\varepsilon}
\newcommand{\glt}{\mathfrak{gl}}
\newcommand{\gt}{\mathfrak{g}}
\newcommand{\htt}{\mathfrak{h}}
\newcommand{\Jo}{\EuScript{J}}
\newcommand{\nt}{\mathfrak{n}}
\newcommand{\Ou}{\mathcal{O}}
\newcommand{\Ro}{\mathcal{R}}
\newcommand{\Rp}{\mathbb{R}}
\newcommand{\ut}{\mathfrak{u}}
\newcommand{\wt}{\widetilde}

\DeclareMathOperator{\aad}{ad}
\DeclareMathOperator{\Ad}{Ad}
\DeclareMathOperator{\GL}{GL}
\DeclareMathOperator{\low}{low}
\DeclareMathOperator{\rk}{rk}
\DeclareMathOperator{\tr}{tr}

\title{%
	Rook placements in $G_2$ and $F_4$ and associated coadjoint orbits
}

\author{%
	Mikhail V. Ignatev,  Matvey A. Surkov
}

\affiliation{
	\address{Mikhail V. Ignatev: Samara National Research University, Ak. Pavlova 1,\break\indent Samara 443011, Russia
	}
	\email{%
		mihail.ignatev@gmail.com
	}
	\address{Matvey A. Surkov: Samara National Research University, Ak. Pavlova 1,\break\indent Samara 443011, Russia
	}
	\email{%
		victorsumaev@yandex.ru
	}
}

\abstract{%
	Let $\mathfrak{n}$ be a maximal nilpotent subalgebra of a simple complex Lie algebra with root system $\Phi$. A subset $D$ of the set $\Phi^+$ of positive roots is called a rook placement if it consists of roots with pairwise non-positive scalar products. To each rook placement $D$ and each map $\xi$ from $D$ to the set $\mathbb{C}^{\times}$ of nonzero complex numbers one can naturally assign the coadjoint orbit $\Omega_{D,\xi}$ in the dual space $\mathfrak{n}^*$. By definition, $\Omega_{D,\xi}$ is the orbit of $f_{D,\xi}$, where $f_{D,\xi}$ is the sum of root covectors $e_{\alpha}^*$ multiplied by $\xi(\alpha)$, $\alpha\in D$ (in fact, almost all coadjoint orbits studied at the moment have such a form, for certain $D$ and $\xi$). It follows from the results of Andr\'e that if $\xi_1$ and $\xi_2$ are distinct maps from $D$ to $\mathbb{C}^{\times}$ then $\Omega_{D,\xi_1}$ and $\Omega_{D,\xi_2}$ do not coincide for classical root systems $\Phi$. We prove that this is true if $\Phi$ is of type $G_2$, or if $\Phi$ is of type $F_4$ and $D$ is orthogonal.
}



\VOLUME{30}
\NUMBER{2}
\YEAR{2022}
\firstpage{129}
\DOI{https://doi.org/10.46298/cm.9041}

\begin{paper}
	
	\section{Introduction and the main result}

	Let $\gt$ be a simple complex Lie algebra, $\bt$ be a Borel subalgebra of $\gt$, $\Phi$ be the root system of $\gt$, $\Phi^+$ be the set of positive roots corresponding to $\bt$, $\nt$ be the nilradical of $\bt$, $N=\exp(\nt)$ be the corresponding nilpotent algebraic group, and $\nt^*$ be the dual space to $\nt$. The group $N$ acts on $\nt$ by the adjoint action; the dual action of $N$ on the space $\nt^*$ is called coadjoint; we will denote the result of this action by $g.\lambda$ for $g\in N$, $\lambda\in\nt^*$. According to the orbit method discovered by A.A. Kirillov in 1962, coadjoint orbits play a key role in representation theory of $N$ (see \cite{Kirillov62}, \cite{Kirillov04}). We will consider a special class of coadjoint orbits defined below.
	
 \begin{definition}
A subset $D$ of $\Phi^+$ is called a \emph{rook placement} if $(\alpha,\beta)\leq0$ for all distinct $\alpha,\beta\in D$, where $(-,-)$ denotes the inner product.
 \end{definition}
	
	The root vectors $e_{\alpha}$, $\alpha\in\Phi^+$ form a basis of $\nt$; we denote by $\{e_{\alpha}^*, \alpha\in D\}$ the dual basis of $\nt^*$. Given a rook placement $D$ and a map $\xi\colon D\to\Cp^{\times}$, we put
	\begin{equation*}
		f_{D,\xi}=\sum_{\alpha\in D}\xi(\alpha)e_{\alpha}^*\in\nt^*.
	\end{equation*}
	
	\begin{definition}
    We say that the coadjoint orbit $\Omega_{D,\xi}$ of the linear form $f_{D,\xi}$ is \emph{associated} with the rook placement $D$ and the map $\xi$.
    \end{definition}
	
	It turns out that almost all coadjoint orbit studied to the moment are associated with certain $D$ and $\xi$ (see, e.g.\ \cite{Andre95}, \cite{AndreNeto06}, \cite{Kostant12}, \cite{Kostant13}, \cite{Panov08}, \cite{IgnatyevPanov09}, \cite{Ignatev09}, \cite{Ignatev11}, \cite{Ignatyev12}). On the other hand, C.A.M. Andr\'e discovered that, for the case of $A_{n-1}$, rook placements themselves provide a nice splitting of $\nt^*$ into a disjoint union of $N$-stable affine subvarieties called basic subvarieties. (We will recall Andr\'e's results in detail in Section~\ref{sect:Andre}, because we will use them for the case of $F_4$.) By definition, the \emph{basic subvariety} $\Ou_{D,\xi}$ corresponding to a rook placement $D$ and a map $\xi\colon D\to\Cp^{\times}$ is $$\Ou_{D,\xi}=\sum_{\alpha\in D}\Omega_{\{\alpha\},\xi_{\alpha}},$$ where $\xi_{\alpha}$ is the restriction of $\xi$ to $\{\alpha\}$.  For $A_{n-1}$, $\nt^*=\bigsqcup_{D,\xi}\Ou_{D,\xi}$ and all $\Ou_{D,\xi}$'s are affine subvarieties of $\nt^*$ (see \cite[Theorem~1]{Andre95}).
	
	Even for $B_n$, $C_n$ and $D_n$, the analogous question is still open. Nevertheless, we may formulate the following conjecture for an arbitrary root system. Non-singularity of a rook placement $D$ means that if $\alpha,\beta\in D$ and $\alpha\neq\beta$ then $\alpha-\beta\notin\Phi^+$; for $A_{n-1}$, all rook placements are automatically non-singular.
	
	\begin{conjecture}Each basic subvariety $\Ou_{D,\xi}$ is an affine subvariety of $\nt^*$\textup, and $$\nt^*=\bigsqcup_{D,\xi}\Ou_{D,\xi},$$ where the union is taken over all non-singular rook placements $D$ and all maps $\xi\colon D\to\Cp^{\times}$.
\end{conjecture}
 
	Direct computations show that this conjecture is true for classical root systems of low rank. In the present paper, we check that this conjecture is true for the case of $G_2$. This is our first main result. In fact, given $D$ and $\xi$, we present an explicit system of equations describing $\Ou_{D,\xi}$.
	
	\begin{theorem}\label{mtheo:G2}
 Let $\Phi=G_2$. Then each basic subvariety $\Ou_{D,\xi}$ is an affine subvariety of $\nt^*$\textup, and $$\nt^*=\bigsqcup_{D,\xi}\Ou_{D,\xi},$$ where the union is taken over all non-singular rook placements $D$ and all maps $\xi\colon D\to\Cp^{\times}$.
 \end{theorem}

	It turns out that, for $A_{n-1}$, if $D$ is a rook placement and $\xi_1$, $\xi_2$ are distinct maps from $D$ to $\Cp^{\times}$ then the associated orbits $\Omega_{D,\xi_1}$ and $\Omega_{D,\xi_2}$ do not coincide (it follows immediately from Andr\'e's theory, since $\Omega_{D,\xi}\subset\Ou_{D,\xi}$, see Section~\ref{sect:Andre}). For other classical root systems this fact can be obtained as a corollary of the case of $A_{n-1}$ (see also \cite{AndreNeto06}). This was used by M.V. Ignatyev and I. Penkov in \cite{IgnatyevPenkov16} and \cite{Ignatyev19} for explicit classification of centrally generated primitive ideals in the universal enveloping algebra $U(\nt)$ for classical root systems.
	
	In \cite{IgnatyevShevchenko21}, M.V. Ignatyev and A.A. Shevchenko, while classifying centrally generated primitive ideals in $U(\nt)$ for exceptional root systems, proved that the analogous is true for certain orthogonal rook placements in $F_4$ and $E_6$, $E_7$, $E_8$. This allows us to formulate the second conjecture for an arbitrary root system.
	
	\begin{conjecture}Let $D$ be a non-singular rook placement and $\xi_1$, $\xi_2$ be distinct maps from $D$ to $\Cp^{\times}$. Then the associated coadjoint orbits $\Omega_{D,\xi_1}$ and $\Omega_{D,\xi_2}$ do not coincide.
 \end{conjecture}

	Our second main result is that this conjecture is true for $F_4$ if $D$ is orthogonal (i.e.\ if all roots from $D$ are pairwise orthogonal).
	
	\begin{theorem}\label{mtheo:F4}
 Let $\Phi=F_4$, $D$ be an orthogonal non-singular rook placement and $\xi_1$, $\xi_2$ be distinct maps from $D$ to $\Cp^{\times}$. Then the associated coadjoint orbits $\Omega_{D,\xi_1}$ and $\Omega_{D,\xi_2}$ do not coincide.
 \end{theorem}

	\section{Andr\'e's theory}\label{sect:Andre}
	
	In this section, we briefly recall Andr\'e's results from \cite{Andre95}, which will be needed in the sequel. Throughout this section, $\Phi$ will be of type $A_{n-1}$. As usual, we identify the set of positive roots with the following subset of the Euclidean space $\Rp^n$:
	$$A_{n-1}^+=\{\epsi_i-\epsi_j, 1\leq i<j\leq n\},$$
	with the standard inner product. Here, $\epsi_1,\ldots,\epsi_n$ denotes the standard basis of $\Rp^n$.
	
	In this case, $\nt$ can be identified with the space of strictly upper-triangular $n\times n$ matrices. Given $\alpha=\epsi_i-\epsi_j\in\Phi^+$, one can pick the $(i,j)$-th  elementary matrix $e_{i,j}$ as a root vector $e_{\alpha}$, so that $[e_{\alpha},e_{\beta}]=\pm e_{\alpha+\beta}$ for $\alpha,\beta\in\Phi^+$ (we put $e_{\alpha+\beta}=0$ if $\alpha+\beta\notin\Phi^+$). We will identify the dual space $\nt^*$ with the space $\nt^-$ of strictly lower-triangular $n\times n$ matrices via the formula $\langle\lambda,x\rangle=\tr(\lambda x)$ for $x\in\nt$, $\lambda\in\nt^-$. The root vectors $e_{\alpha}$, $\alpha\in\Phi^+$ form a basis of $\nt$; let $\{e_{\alpha}^*, \alpha\in\Phi^+\}$ be the dual basis of $\nt^*$ (in fact, $e_{i,j}^*=e_{j,i}$).
	
	The group $N$ is the group of all upper-triangular $n\times n$ matrices with $1$'s on the diagonal. It acts on its Lie algebra $\nt$ via the adjoint action $\Ad_g(x)=gxg^{-1}$, $g\in N$, $x\in\nt$. The dual action of $N$ on the space $\nt^*$ is called coadjoint; we will denote the result of this action by $g.\lambda$, $g\in N$, $\lambda\in\nt^*$. It is easy to check that, after the identification of $\nt^*$ with $\nt^-$, this action has the form $g.\lambda=(g\lambda g^{-1})_{\low}$. Here, given an $n\times n$ matrix $a$, we set
	\begin{equation*}
		(a_{\low})_{i,j}=\begin{cases} a_{i,j},&\text{if }i>j,\\
			0&\text{otherwise}.
		\end{cases}
	\end{equation*}
	
	\begin{definition}Pick a number $k$ from $1$ to $n$. We call the sets
		\begin{equation*}
			\Ro_k=\{\epsi_j-\epsi_k, 1\leq j<k\}, \Co_k=\{\epsi_k-\epsi_i, k<i\leq n\}
		\end{equation*}
		the \emph{$k$-th row} and the \emph{$k$-th column} of $\Phi^+$, respectively. We say that the number $i$ (respectively, the number $j$) is the \emph{row} (respectively, the \emph{column}) of a root $\alpha=\epsi_i-\epsi_j$.
  \end{definition}

    \begin{example}
    Let $n=6$. On the picture below boxes from $\Ro_5\cup\Co_2$ are grey. Here we identify a root $\epsi_i-\epsi_j\in\Phi^+$ with the box $(j,i)$.
		\begin{center}
			\includegraphics{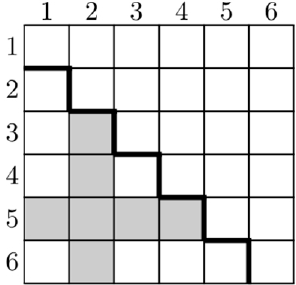}
		\end{center}
	\end{example}
	
	To each rook placement $D\subset\Phi^+$ and each map $\xi\colon D\to\Cp^{\times}$, one can assign the linear form $$f_{D,\xi}=\sum_{\alpha\in D}\xi(\alpha)e_{\alpha}^*\in\nt^*.$$
	
	\begin{example}Let $n=8$, $D=\{\epsi_1-\epsi_3, \epsi_2-\epsi_6, \epsi_3-\epsi_7, \epsi_4-\epsi_5, \epsi_6-\epsi_8\}$. On the picture below we schematically draw the linear form $f_{D,\xi}$ putting symbols $\otimes$ in the boxes corresponding to the roots from $D$.
		\begin{center}
			\includegraphics{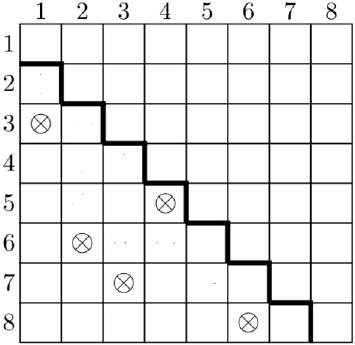}
		\end{center}
	\end{example}
	
	It follows immediately from the definition of a rook placement that $$|D\cap\Ro_k|\leq1, |D\cap\Co_k|\leq1$$ for all $k$. This explains the term ``rook placement'': if we identify symbols $\otimes$ from $f_{D,\xi}$ with rooks on the lower-triangular chessboard, then these rooks do not hit each other.
	
	Now, given a root $\alpha\in D$, we denote by $\xi_{\alpha}$ the restriction of the map $\xi$ to the subset $\{\alpha\}$, and put
	\begin{equation*}
		\Ou_{D,\xi}=\sum_{\alpha\in D}\Ou_{\{\alpha\},\xi_{\alpha}}.
	\end{equation*}
	Clearly, $\Omega_{D,\xi}\subset\Ou_{D,\xi}$.
	
	\begin{definition}The set $\Ou_{D,\xi}$ is called a \emph{basic subvariety} of $\nt^*$ corresponding to the rook placement $D$ and the map $\xi$.
 \end{definition}
	
	Accordingly to \cite[Theorem~1]{Andre95}, $\nt^*$ is a disjoint unions of basic subvarieties:
	\begin{equation*}
		\nt^*=\bigsqcup_{D,\xi}\Ou_{D,\xi},
	\end{equation*}
	where the union is taken over all rook placements in $A_{n-1}^+$ and all maps $\xi\colon D\to\Cp^{\times}$. Formally, Andr\'e considered the case of finite ground field, but his proofs are valid over $\Cp$, too. Furthermore, each basic subvariety $\Ou_{D,\xi}$ is in fact an affine subvariety of $\nt^*$, and Andr\'e presented an explicit set of defining equations for it. To describe this set, we need some more notation.
	
	\begin{definition}A root $\alpha\in\Phi^+$ is called $\beta$-\emph{singular} for a root $\beta\in\Phi^+$ if $\beta-\alpha\in\Phi^+$. The set of all $\beta$-singular roots is denoted by $S(\beta)$.
 \end{definition}

	Let $D$ be a rook placement and $\xi\colon D\to\Cp^{\times}$ be a map. We denote $$S(D)=\bigcup_{\beta\in D}S(\beta)$$ and $R(D)=\Phi^+\setminus S(D)$. Obviously, $D\subset R(D)$. Roots from $S(D)$ (respectively, from $R(D)$) are called $D$-\emph{singular} (respectively, $D$-\emph{regular}).
	
	\begin{example}Let $n=10$, $D=\{\epsi_1-\epsi_6, \epsi_3-\epsi_{10}, \epsi_5-\epsi_8\}$. On the picture below, boxes corresponding to the roots from $D$ are filled by $\otimes$'s, as above, while the boxes corresponding by the $D$-singular roots are marked gray.\end{example}
	
	\begin{center}
		\includegraphics{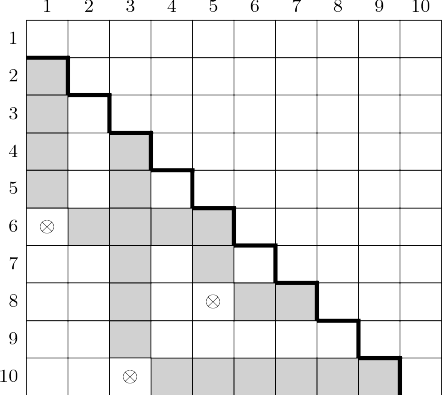}
	\end{center}
	
	It turns out that to each $D$-regular root $\alpha$ one can assign the defining equation of $\Ou_{D,\xi}$. Namely, there is a natural partial order on $\Phi^+$: we write $\alpha\leq\beta$ if $\beta-\alpha$ is a sum of positive roots. Evidently, $\epsi_i-\epsi_j\leq\epsi_r-\epsi_s$ if and only if $s\leq j$ and $r\geq i$ (in other words, on our pictures $\epsi_r-\epsi_s$ is located non-strictly to the South-West from $\epsi_i-\epsi_j$). Given $\alpha\in R(D)$, we set $D(\alpha)=\{\alpha\}\cup\{\beta\in D\mid\beta\geq\alpha\}.$ Now, let $R_D(\alpha)$ (respectively $C_D(\alpha)$) be the set of all rows (respectively, of all columns) of the roots from $D(\alpha)$. Finally, for a matrix $\lambda\in\nt^*$, we denote by $\Delta^D_{\alpha}(\lambda)$ the minor of the matrix $\lambda$ with the set of rows $R_D(\alpha)$ and the set of columns $C_D(\alpha)$. We assume that the numbers of rows and columns are taken in the increasing order.
	
	For instance, in the previous example, if $\alpha=\epsi_4-\epsi_5$ then
    \[
        D(\alpha)
        = \{\epsi_1-\epsi_6, \epsi_3-\epsi_{10}, \epsi_4-\epsi_5\},
    \]
    hence $R_D(\alpha)=\{5, 6, 10\}$, $C_D(\alpha)=\{1, 3, 4\}$ and
	\begin{equation*}
		\Delta^D_{\alpha}(\lambda)=\begin{vmatrix}
			\lambda_{5,1}&\lambda_{5,3}&\lambda_{5,4}\\
			\lambda_{6,1}&\lambda_{6,3}&\lambda_{6,4}\\
			\lambda_{10,1}&\lambda_{10,3}&\lambda_{10,4}\\
		\end{vmatrix}.
	\end{equation*}
	On the other hand, for $\alpha=\epsi_5-\epsi_8\in D$, one has $D(\alpha)=\{\epsi_3-\epsi_{10}, \epsi_5-\epsi_8\}$, $R_D(\alpha)=\{8, 10\}$, $C_D(\alpha)=\{3, 5\}$, and, consequently,
	\begin{equation*}
		\Delta^D_{\alpha}(\lambda)=\begin{vmatrix}
			\lambda_{8,3}&\lambda_{8,5}\\
			\lambda_{10,3}&\lambda_{10,5}\\
		\end{vmatrix}.
	\end{equation*}
	
	Thanks to \cite[Proposition~2]{Andre95}, a matrix $\lambda\in\nt^*$ belongs to $\Ou_{D,\xi}$ if and only if
	\begin{equation*}
		\Delta^D_{\alpha}(\lambda)=\Delta^D_{\alpha}(f_{D,\xi})\text{ for all }\lambda\in R(D).
	\end{equation*}
	Precisely, $\Delta^D_{\alpha}(\lambda)=0$ for all $\alpha\in R(D)\setminus D$ and $\Delta^D_{\alpha}(\lambda)=\pm\prod_{\beta\in D(\alpha)}\xi(\beta)$ for $\alpha\in D$. It follows immediately that $$\dim\Ou_{D,\xi}=|S(D)|.$$

 \begin{remark}\label{nota:Andre}
 Actually, Andr\'e's proof of the fact that each $\lambda\in\nt^*$ belongs to a certain basic subvariety $\Ou_{D,\xi}$ is very straightforward. Namely, there is a total order $\leq_t$ on $\Phi^+$ refining the partial order $\leq$ defined above. By definition, $\epsi_i-\epsi_j<_t\epsi_r-\epsi_s$ if $s<j$ or $s=j$, $i<r$. Now, given $\lambda\in\nt^*$, we inductively construct $\Ou_{D,\xi}$ containing $\lambda$. If $\lambda=0$, then $D=\varnothing$ and $\xi$ is the unique empty map from $\varnothing$ to $\Cp^{\times}$. If $\lambda\neq0$ then we find the smallest (with respect to $\leq_t$) root $\alpha_1$ from $\Phi^+$ such that $\lambda(e_{\alpha_1})\neq0$, and put $D=\{\alpha_1\}$, $\xi(\alpha_1)=\lambda(\alpha_1)$. If $\lambda\in\Ou_{D,\xi}$, then we are done. Otherwise, let $\alpha_2$ be the smallest root such that $\Delta^D_{\alpha_2}(\lambda)\neq\Delta^D_{\alpha_2}(f_{D,\xi})$. Then we add $\alpha_2$ to $D$ and define $\xi(\alpha_2)$ in the obvious way. Now, one can repeat this procedure to obtain the required basic subvariety $\Ou_{D,\xi}$.
 \end{remark}
	
	\section{Case $\Phi=G_2$}
	
	In this section, we prove our first main result, Theorem~\ref{mtheo:G2}. First, we briefly recall some basic facts about the simple Lie algebra $\gt$ of type $G_2$ and its maximal nilpotent subalgebra $\nt$. By definition, the root system $\Phi=G_2$ has the form $\Phi=\Phi^+\cup\Phi^-$, where $\Phi^+=\{ \alpha, \beta, \alpha+\beta, 2\alpha+\beta, 3\alpha+\beta, 3\alpha+2\beta \}$, $\Phi^-=-\Phi^+$,
	and $\alpha$, $\beta$ are vectors from $\Rp^2$ such that $||\alpha||^2=1, ||\beta||^2=3$ and the angle between $\alpha$ and $\beta$ equals $5\pi/6$. There is a Cartan decomposition $\gt=\htt\oplus\nt\oplus\nt_-$, where $\htt$ is a Cartan subalgebra of $\gt$, and $\nt$ has a basis consisting of the root vectors $e_{\gamma}$, $\gamma\in\Phi^+$.
	
	It is well known that there exists nonzero scalars $c_i$, $1\leq i\leq5$, such that
 \begin{align*}
     &[e_\alpha, e_\beta] = c_1\cdot e_{\alpha+\beta},\\
     &[e_\alpha, e_{\alpha+\beta}] = c_2\cdot e_{2\alpha+\beta} \\
     &[e_\alpha, e_{2\alpha+\beta}] = c_3\cdot e_{3\alpha+\beta}\\
     &[e_{3\alpha+\beta}, e_\beta] = c_4\cdot e_{3\alpha+2\beta} \\
     &[e_{\alpha+\beta}, e_{2\alpha+\beta}] = c_5\cdot e_{3\alpha+2\beta}.
 \end{align*}
	In fact, one can choose the root vectors so that $c_1=1$, $c_2=2$, $c_3=3$, $c_4=1$, $c_5=3$, but we will not use these explicit values in the sequel. One can immediately check that $c_1c_5=c_3c_4$ for an arbitrary choice of the root vectors.
	
	Recall the definition of the group $N=\exp(\nt)$ and the coadjoint action of $N$ on the dual space $\nt^*$. It is straightforward to check that this action has the form
	\begin{align*}
	    (\exp(x).\lambda)(y)&=\lambda(\exp(-\ad{x})(y)) \\
     &=\lambda(y) - \lambda([x, y]) + \dfrac{1}{2!}\lambda([x, [x, y]]) - \ldots,
	\end{align*}
for $ x,y\in\nt, \lambda\in\nt^*$.
	
	Now, let $D$ be a non-singular rook placement in $\Phi^+$. Recall that non-singularity means that $\gamma\notin S(\delta)$ for all distinct $\gamma, \delta\in D$, where $S(\delta)$ denotes the set of $\delta$-singular roots in $\Phi^+$.
	Fix a map $\xi\colon D\to\Cp^{\times}$, and recall that, by definition, $\Omega_{D,\xi}$ is the coadjoint orbit of the linear form $f_{D,\xi}$. It follows immediately that if $\gamma$ is a maximal (with respect to the partial order $\leq$ on $\Phi^+$) among all roots from $D$ then $\lambda(e_{\gamma})=f_{D,\xi}(e_{\gamma})=\xi(\gamma)
	$
	for all $\lambda\in\Omega_{D,\xi}$. Similarly, $\lambda(e_{\gamma})=0$ for all $\lambda\in\Omega_{D,\xi}$, if there are no $\delta\in D$ such that $\delta\geq\gamma$. Recall also the definition of $\Ou_{D,\xi}$.
	
	Given $\gamma\in\Phi^+$, we write $\lambda_{\gamma}=\lambda(e_{\gamma})$, so that $\lambda=\sum_{\gamma\in\Phi^+}\lambda_{\gamma}e_{\gamma}^*$. We will prove Theorem\ref{mtheo:G2} as an immediate corollary of the following key proposition:

	\begin{proposition}\label{prop:G2_O_D_xi}Let $D$ be a non-singular rook placement in $\Phi^+$\textup, $\xi\colon D\to\Cp^\times$ be a map. Pick a linear form $\lambda\in\nt^*$. Then $\lambda\in\Ou_{D,\xi}$ if and only if $\lambda$ satisfy the following system of equations.
		{
\begin{longtable}{||l|l|l||}
			\hline
			&\begin{tabular}[t]{l}$D$\end{tabular}&\begin{tabular}[t]{l}\textup{System of equations for} $\Ou_{D,\xi}$\end{tabular}
			\\
            \hline\hline
			\textup{1}
			&$\varnothing$
			&
			\hspace{-0.3cm}\begin{tabular}[t]{l}
				$\lambda_{\gamma}=0$ for all $\gamma\in\Phi^+$
			\end{tabular} \\
   \hline
			\textup{2}
			&$\alpha$
			&
			\hspace{-0.3cm}\begin{tabular}[t]{l}
				$\lambda_{\alpha}=\xi(\alpha)$,\\
				$\lambda_{\gamma}=0$ for $\gamma\neq\alpha$
			\end{tabular}
			\\
			\hline
			\textup{3}
			&$\beta$
			&
			\hspace{-0.3cm}\begin{tabular}[t]{l}
				$\lambda_{\beta}=\xi(\beta)$,\\
				$\lambda_{\gamma}=0$ for $\gamma\neq\beta$
			\end{tabular}\\
   \hline
			\textup{4}
			&$\alpha+\beta$
			&
			\hspace{-0.3cm}\begin{tabular}[t]{l}
				$\lambda_{\alpha+\beta}=\xi(\alpha+\beta)$,\\
				$\lambda_{2\alpha+\beta}=\lambda_{3\alpha+\beta}=\lambda_{3\alpha+2\beta}=0$ 
			\end{tabular}
			\\
			\hline
			\textup{5}
			&$2\alpha+\beta$
			&
			\hspace{-0.3cm}\begin{tabular}[t]{l}
				$\lambda_{2\alpha+\beta}=\xi(2\alpha+\beta)$,\\
				$2c_2\lambda_{\beta}\lambda_{2\alpha+\beta}-c_1\lambda_{\alpha+\beta}^2=0$,\\
				$\lambda_{3\alpha+\beta}=\lambda_{3\alpha+2\beta}=0$
			\end{tabular}
			\\
   \hline
   \textup{6}
			&$3\alpha+\beta$
			&
			\hspace{-0.3cm}\begin{tabular}[t]{l}
				$6c_3^2\lambda_{\beta}\lambda_{3\alpha+\beta}^2-c_1c_2\lambda_{2\alpha+\beta}^3=0$,\\
				$2c_3\lambda_{\alpha+\beta}\lambda_{3\alpha+\beta}-c_2\lambda_{2\alpha+\beta}^2=0$,\\
				$\lambda_{3\alpha+\beta}=\xi(3\alpha+\beta)$,\\ $\lambda_{3\alpha+2\beta}=0$
			\end{tabular}
			\\
			\hline
			\textup{7}
			&$3\alpha+2\beta$
			&
			\hspace{-0.3cm}\begin{tabular}[t]{l}
				$2c_5\lambda_{\alpha}\lambda_{3\alpha+2\beta}-2c_3\lambda_{\alpha+\beta}\lambda_{3\alpha+\beta}+c_2\lambda_{2\alpha+\beta}^2=0$,\\
				$\lambda_{3\alpha+2\beta}=\xi(3\alpha+2\beta)$
			\end{tabular}
			\\
   \hline
   \textup{8}
			&$\alpha, \beta$
			&
			\hspace{-0.3cm}\begin{tabular}[t]{l}
				$\lambda_{\alpha}=\xi(\alpha)$, \\ $\lambda_{\beta}=\xi(\beta)$,\\
				$\lambda_{\gamma}=0$ for $\gamma\neq\alpha, \beta$
			\end{tabular}
			\\
			\hline
			\textup{9}
			&$\beta, 2\alpha+\beta$
			&
			\hspace{-0.3cm}\begin{tabular}[t]{l}
				$\lambda_{2\alpha+\beta}=\xi(2\alpha+\beta)$,\\
				$2c_2\lambda_{\beta}\lambda_{2\alpha+\beta}-c_1\lambda_{\alpha+\beta}^2=$\\
				$2c_2\xi(\beta)\xi(2\alpha+\beta)$,\\
				$\lambda_{3\alpha+\beta}=\lambda_{3\alpha+2\beta}=0$
			\end{tabular}\\
            \hline
			\textup{10}
			&$\beta, 3\alpha+\beta$
			&
			\hspace{-0.3cm}\begin{tabular}[t]{l}
				$6c_3^2\lambda_{\beta}\lambda_{3\alpha+\beta}^2-c_1c_2\lambda_{2\alpha+\beta}^3=6c_3^2\xi(\beta)\xi(3\alpha+\beta)^2$,\\
				$2c_3\lambda_{\alpha+\beta}\lambda_{3\alpha+\beta}-c_2\lambda_{2\alpha+\beta}^2=0$,\\
				$\lambda_{3\alpha+\beta}=\xi(3\alpha+\beta)$,\\
				$\lambda_{3\alpha+2\beta}=0$
			\end{tabular}
			\\
			\hline
			\textup{11}
			&$\alpha+\beta, 3\alpha+\beta$
			&
			\hspace{-0.3cm}
   \begin{tabular}[t]{l}
				$2c_3\lambda_{\alpha+\beta}\lambda_{3\alpha+\beta}-c_2\lambda_{2\alpha+\beta}^2=2c_3\xi(\alpha+\beta)\xi(3\alpha+\beta)$,\\
				$\lambda_{3\alpha+\beta}=\xi(3\alpha+\beta)$,\\
				$\lambda_{3\alpha+2\beta}=0$
			\end{tabular}\\
			\hline
            \textup{12}
			&$\alpha, 3\alpha+2\beta$
			&
			\hspace{-0.3cm}\begin{tabular}[t]{l}
				$2c_5\lambda_{\alpha}\lambda_{3\alpha+2\beta}-2c_3\lambda_{\alpha+\beta}\lambda_{3\alpha+\beta}+c_2\lambda_{2\alpha+\beta}^2=c_5\xi(\alpha)\xi(3\alpha+2\beta)$,\\
				$\lambda_{3\alpha+2\beta}=\xi(3\alpha+2\beta)$
			\end{tabular}
			\\
			\hline
	\end{longtable}}\end{proposition}
	\begin{proof}The proof will be performed for all rook placements in $\Phi^+$ subsequently. First, assume that $|D|=1$, and in that case, $\Omega_{D,\xi}=\Ou_{D,\xi}$. Pick a linear form $\lambda\in\Ou_{D,\xi}$. Then there exists $x=\sum_{\gamma\in\Phi^+}x_{\gamma}e_{\gamma}\in\nt$ such that $\lambda=\exp(x).f_{D,\xi}$. Let us proceed case-by-case. Cases 1, 2, 3 from the table above are evident, so we start from case 4.
		
		Case 4: $D = \{ \alpha + \beta \}$.   
  
  It follows immediately from the paragraph before Proposition~\ref{prop:G2_O_D_xi} that $\lambda_{\alpha+\beta} = \xi(\alpha + \beta)$ and $\lambda_{2\alpha+\beta} = \lambda_{3\alpha+\beta} = \lambda_{3\alpha+2\beta} = 0$. To compute $\lambda_{\alpha}$, we note that, obviously, $(\alpha+\beta)-\alpha$ can be uniquely represented as a sum of positive roots: $(\alpha+\beta)-\alpha=\beta$. Therefore,
		\begin{equation*}
			\lambda_{\alpha} = \lambda(e_\alpha) = \xi(\alpha + \beta)e_{\alpha + \beta}^*(e_\alpha - [x, e_\alpha]) = \xi(\alpha + \beta)e_{\alpha + \beta}^*(-[x_\beta e_\beta, e_\alpha]) = \xi(\alpha + \beta)c_1x_\beta.
		\end{equation*}
		Similarly, $\lambda_{\beta} = -\xi(\alpha + \beta)c_1x_\alpha.$ Since $x_{\alpha}$ and $x_{\beta}$ can be arbitrary, we obtain the required system of equations.
		
		Case 5: $D = \{ 2\alpha + \beta \}$.
  
  Here $\lambda_{2\alpha+\beta} = \xi(2\alpha + \beta)$ and $ \lambda_{3\alpha+\beta} = \lambda_{3\alpha+2\beta} = 0$. Now, $(2\alpha+\beta)-(\alpha+\beta)=\alpha$ is the unique representation of $(2\alpha+\beta)-(\alpha+\beta)$ as a sum of positive roots, hence
		\begin{equation*}
			\lambda_{\alpha+\beta} = \lambda(e_{\alpha + \beta}) = \xi(2\alpha + \beta)e_{2\alpha + \beta}^*(e_{\alpha + \beta} - [x, e_{\alpha + \beta}]) = -\xi(2\alpha + \beta)c_2x_\alpha.
		\end{equation*}
		Next, since $(2\alpha + \beta)-\beta=2\alpha$ is the unique representation of $(2\alpha+\beta)-\beta$ as a sum of positive roots, we obtain
		\begingroup
    \begin{align*}
				\lambda_{\beta}&= \lambda(e_\beta) \\
    &= \xi(2\alpha + \beta)e_{2\alpha + \beta}^*(e_\beta - [x, e_\beta] + \frac{1}{2!}[x, [x, e_\beta]]) \xi(2\alpha + \beta)e_{2\alpha + \beta}^*(\frac{1}{2}[x_\alpha e_\alpha, [x_\alpha e_\alpha, e_\beta]])\\
				&= \frac{1}{2}\xi(2\alpha + \beta)c_1x_\alpha^2 e_{2\alpha + \beta}^*([e_\alpha, e_{\alpha + \beta}]) = \frac{1}{2}\xi(2\alpha + \beta)c_1c_2x_\alpha^2.
		\end{align*}
    \endgroup
		Finally, since $2\alpha\notin\Phi^+$, we can obtain $2\alpha+\beta$ either by adding to $\alpha$ the roots $\beta$ and $\alpha$ subsequently, or by adding to $\alpha$ the root $\alpha+\beta$. So,
		\begin{equation*}
			\lambda_{\alpha} = \lambda(e_\alpha) = \xi(2\alpha + \beta)e_{2\alpha + \beta}^*(e_\alpha - [x, e_\alpha] + \frac{1}{2!}[x, [x, e_\alpha]]) =\xi(2\alpha + \beta) (c_2x_{\alpha + \beta} - \frac{1}{2}c_1c_2x_\alpha x_\beta).
		\end{equation*}
		Thus, $\lambda_{\alpha}$ can be arbitrary, while $2c_2\lambda_{2\alpha+\beta}\lambda_{\beta}=c_1\lambda_{\alpha+\beta}^2$, as required.
		
		Case  6: $D = \{ 3\alpha + \beta \}$.
  
  Arguing as above, we see that $\lambda_{3\alpha+\beta} = \xi(3\alpha + \beta)$, $\lambda_{3\alpha+2\beta} = 0$,
		\begin{align*}
				&\lambda_{2\alpha+\beta} = -\xi(3\alpha + \beta)c_3x_\alpha, \lambda_{\alpha+\beta}= \frac{1}{2}\xi(3\alpha + \beta)c_2c_3x_\alpha^2, \lambda_{\beta}=-\frac{1}{6}\xi(3\alpha + \beta)c_1c_2c_3x_\alpha^3,\\
				&\lambda_{\alpha}=\xi(3\alpha + \beta) (c_3x_{2\alpha + \beta} - \frac{1}{2}c_2c_3x_\alpha x_{\alpha + \beta} + \frac{1}{6}c_1c_2c_3x_\alpha^2x_\beta).
		\end{align*}
		Now, it is clear that the equations from the table above define $\Ou_{D,\xi}$.

		Case 7: $D = \{ 3\alpha + 2\beta \}$. 
  
  Here $\lambda_{3\alpha+2\beta} = \xi(3\alpha + 2\beta)$, $\lambda_{3\alpha+\beta} = \xi(3\alpha + 2\beta)c_4x_\beta,$
		\begin{align*}
				&\lambda_{2\alpha+\beta} = \xi(3\alpha + 2\beta) (-c_5x_{\alpha + \beta} - \frac{1}{2}c_3c_4x_\alpha x_\beta),\\
				&\lambda_{\alpha+\beta} = \xi(3\alpha + 2\beta) (c_5x_{2\alpha + \beta} + \frac{1}{2}c_2c_5x_\alpha x_{\alpha + \beta} + \frac{1}{6}c_2c_3c_4x_\alpha^2x_\beta),\\
				&\lambda_{\beta} = \xi(3\alpha + 2\beta)(-c_4x_{3\alpha + \beta} - \frac{1}{2}c_1c_5x_\alpha x_{2\alpha + \beta} - \frac{1}{6}c_1c_2c_5x_\alpha^2 x_{\alpha + \beta} - \frac{1}{24}c_1c_2c_3c_4x_\alpha^3 x_\beta),\\
				&\lambda_{\alpha} = \xi(3\alpha + 2\beta)(c_1c_5x_\beta x_{2\alpha + \beta} - \frac{1}{2}c_2c_5x_{\alpha + \beta}^2 + \frac{1}{24}c_1c_2c_3c_4x_\alpha^2 x_\beta^2).
		\end{align*}
		One can immediately check that $\lambda$ satisfies the required system of equations. On the other hand, given arbitrary $\lambda_{2\alpha+\beta}$, $\lambda_{\alpha+\beta}$, $\lambda_{\beta}$, one can put $x_\beta= \dfrac{\lambda_{3\alpha+\beta}}{c_4\xi(3\alpha+2\beta)}$,
		\begin{align*}
				x_{\alpha+\beta}&= -\frac{\lambda_{2\alpha+\beta} + \frac{1}{2}c_3c_4x_\alpha x_\beta\xi(3\alpha+2\beta)}{c_5\xi(3\alpha+2\beta)}=-\frac{\lambda_{2\alpha+\beta} + \frac{1}{2}c_3\lambda_{3\alpha+\beta}x_\alpha}{c_5\xi(3\alpha+2\beta)},\\
				x_{2\alpha+\beta}&= \frac{\lambda_{\alpha+\beta} - \frac{1}{2}c_2c_5x_\alpha x_{\alpha+\beta}\xi(3\alpha+2\beta) - \frac{1}{6}c_2c_3c_4x_\alpha^2 x_\beta\xi(3\alpha+2\beta)}{c_5\xi(3\alpha+2\beta)}\\
				&=\frac{\lambda_{\alpha+\beta} + \frac{1}{2}c_2x_\alpha(\lambda_{2\alpha+\beta}+\frac{1}{2}c_3\lambda_{3\alpha+\beta}x_\alpha)  - \frac{1}{6}c_2c_3\lambda_{3\alpha+\beta}x_\alpha^2}{c_5\xi(3\alpha+2\beta)}\\
				&=
				\frac{\lambda_{\alpha+\beta} + \frac{1}{2}c_2x_\alpha\lambda_{2\alpha+\beta} + \frac{1}{12}c_2c_3\lambda_{3\alpha+\beta}x_\alpha^2}{c_5\xi(3\alpha+2\beta)}.
		\end{align*}
		It is straightforward to check that, for these values of $x_{\beta}$, $x_{\alpha+\beta}$ and $x_{2\alpha+\beta}$, one has $$\lambda_{\alpha}=\dfrac{c_3\lambda_{\alpha+\beta}\lambda_{3\alpha+\beta}}{c_5\lambda_{3\alpha+2\beta}} - \dfrac{c_2\lambda_{2\alpha+\beta}^2}{2c_5\lambda_{3\alpha+2\beta}}.$$
		Thus, $\Ou_{D,\xi}=\Omega_{D,\xi}$ is exactly the set of solutions of the required system of equations.
		
		Cases 8--12 can be considered uniformly (in all these cases $|D|=2$). In all cases, except case 11, $D$ contains a basis root $\gamma$ ($\gamma=\alpha$ or $\gamma=\beta$). Since the coadjoint orbit of $\xi(\gamma)e_{\gamma}^*$ is $\{\xi(\gamma)e_{\gamma}^*\}$, everything is evident.  Case 11 is an easy exercise.\end{proof}
	
	We are now ready to prove our first main result, Theorem~\ref{mtheo:G2}, which claims that, for $\Phi=G_2$,  $\nt^*=\bigsqcup_{D,\xi}\Ou_{D,\xi}$, where the union is taken over all non-singular rook placements $D$ and all maps $\xi\colon D\to\Cp^{\times}$.
	
	\begin{proof}[Proof of Theorem~\ref{mtheo:G2}] Using Proposition~\ref{prop:G2_O_D_xi}, one can check that each $\lambda\in\nt^*$ belongs to exactly one orbit $\Ou_{D,\xi}$. Namely,  pick a linear form $\lambda \in \mathfrak{n}^*$. Then exactly one of the following cases can occur.
	
	\begin{itemize}
		\item $\lambda_{\alpha}=\lambda_{\beta}=\lambda_{\alpha+\beta}=\lambda_{2\alpha+\beta}=\lambda_{3\alpha+\beta}=\lambda_{3\alpha+2\beta}=0$. Then $\lambda\in\Ou_{\varnothing,\xi}$.
		
		\item $\lambda_{\beta}=\lambda_{\alpha+\beta}=\lambda_{2\alpha+\beta}=\lambda_{3\alpha+\beta}=\lambda_{3\alpha+2\beta}=0, \lambda_{\alpha}\neq 0$. Then $\lambda\in\Ou_{\{ \alpha \},\xi}$.
		
		\item $\lambda_{\alpha+\beta}=\lambda_{2\alpha+\beta}=\lambda_{3\alpha+\beta}=\lambda_{3\alpha+2\beta}=0, \lambda_{\beta}\neq 0$. If $\lambda_{\alpha} = 0$ then $\lambda\in\Ou_{\{ \beta \},\xi}$, if $\lambda_{\alpha} \neq 0$ then $\lambda\in\Ou_{\{ \alpha, \beta \},\xi}$.
		
		\item $\lambda_{2\alpha+\beta}=\lambda_{3\alpha+\beta}=\lambda_{3\alpha+2\beta}=0, \lambda_{\alpha+\beta}\neq 0$. Then $\lambda\in\Ou_{\{ \alpha+\beta \},\xi}$. 
		
		\item $\lambda_{3\alpha+\beta}=\lambda_{3\alpha+2\beta}=0, \lambda_{2\alpha+\beta}\neq 0$. If $\lambda_{\beta} = \dfrac{c_1\lambda_{\alpha+\beta}^2}{2c_2\lambda_{2\alpha+\beta}}$ then $\lambda\in\Ou_{\{ 2\alpha + \beta \},\xi}$, otherwise $\lambda\in\Ou_{\{ \beta,2\alpha+\beta \},\xi}$.
		
		\item $\lambda_{3\alpha+2\beta}=0, \lambda_{3\alpha+\beta}\neq 0$. If $\lambda_{\beta} = \dfrac{c_1c_2\lambda_{2\alpha+\beta}^2}{6c_3^2\lambda_{3\alpha+\beta}^2}$ and $\lambda_{\alpha+\beta} = \dfrac{c_2\lambda_{2\alpha+\beta}^2}{2c_3\lambda_{3\alpha+\beta}}$ then $\lambda\in\Ou_{\{ 3\alpha+\beta \},\xi}$. If $\lambda_{\beta} \neq \dfrac{c_1c_2\lambda_{2\alpha+\beta}^2}{6c_3^2\lambda_{3\alpha+\beta}^2}$ and $\lambda_{\alpha+\beta} = \dfrac{c_2\lambda_{2\alpha+\beta}^2}{2c_3\lambda_{3\alpha+\beta}}$ then $\lambda\in\Ou_{\{ \beta,3\alpha+\beta \},\xi}$. If $\lambda_{\alpha+\beta} \neq \dfrac{c_2\lambda_{2\alpha+\beta}^2}{2c_3\lambda_{3\alpha+\beta}}$ then $\lambda\in\Ou_{\{ \alpha+\beta,3\alpha+\beta \},\xi}$.
		
		\item $\lambda_{3\alpha+2\beta} \neq 0$. If $\lambda_{\alpha} = \dfrac{c_3\lambda_{\alpha+\beta}\lambda_{3\alpha+\beta}}{c_5\lambda_{3\alpha+2\beta}} - \dfrac{c_2\lambda_{2\alpha+\beta}^2}{2c_5\lambda_{3\alpha+2\beta}}$ then $\lambda\in\Ou_{\{ 3\alpha+2\beta \},\xi}$, otherwise $\lambda\in\Ou_{\{ \alpha, 3\alpha+2\beta \},\xi}$.
	\end{itemize}
	
	The proof is complete.
 \end{proof}
	
	\begin{remark}
 \begin{enumerate}[i)]
     \item In fact, $\Ou_{D,\xi}=\Omega_{D,\xi}$ for all $D$ (and $\xi$), except case 11.

     \item There exists exactly one singular rook placement in $G_2^+$, namely, $D=\{\alpha, \alpha+\beta\}$. We do not consider this rook placement because
     \[
     \Omega_{D,\xi}
     =\Ou_{D,\xi}
     =\Omega_{\{\alpha+\beta\},\xi_{\alpha+\beta}}
     =\Ou_{\{\alpha+\beta\},\xi_{\alpha+\beta}},
     \]
     where $\xi_{\alpha+\beta}=\restr{\xi}{\{\alpha+\beta\}}$.

     \item It follows from Proposition~\ref{prop:G2_O_D_xi} that $\dim\Ou_{D,\xi}=|S(D)|$ does not depend on $\xi$, as for $A_{n-1}^+$.
 \end{enumerate}
 \end{remark}

	\section{Case $\Phi=F_4$}
	
	In this section we prove our second main result, Theorem~\ref{mtheo:F4}. To do this, we firstly prove the following simple lemma. Let $D$ be a non-singular orthogonal rook placement in $\Phi^+$, where $\Phi=F_4$, and $\xi_1, \xi_2\colon D\to\Cp^{\times}$ be a map. Assume that there is the unique maximal root $\beta_0$ in $D$ (with respect to the natural order on $\Phi^+$).
	
	\begin{lemma}\label{lemm:F4_submax}
 Let $\beta\in D\setminus\{\beta_0\}$ be such that $\gamma\ngtr\beta$ for all $\gamma\in D\setminus\{\beta_0\}$. Assume that $\beta_0-\beta=\gamma_1+\ldots+\gamma_k$ can be uniquely expressed as a sum of positive roots $\gamma_j$\textup, $1\leq j\leq k$. Further\textup, assume that $$\beta+\sum_{j\in J}\gamma_j\in\Phi^+$$ for each subset $J\subset\{1, \ldots, k\}$. If $\xi_1(\beta)\neq\xi_2(\beta)$, then $\Omega_{D,\xi_1}$ and $\Omega_{D,\xi_2}$ do not coincide.\end{lemma}
 
	\begin{proof}Suppose that $\Omega_{D,\xi_1} = \Omega_{D,\xi_2}$. Then $\xi_1(\beta_0)=\xi_2(\beta_0)$. Let $\xi\colon D\to\Cp^{\times}$ be a map. Pick an element $x\in\nt$ and denote $\mu=\exp{x}.f_{D,\xi}$. One has
			\begin{align*}
				\mu(e_{\beta_0-\gamma_j})&=(\exp x.f_{D,\xi})(e_{\beta_0-\gamma_j})\\
                & = f_{D,\xi}(e_{\beta_0-\gamma_j}-[x,e_{\beta_0-\gamma_j}]+\ldots)\\
				&=a-f_{D,\xi}([x_{\gamma_j}e_{\gamma_j},e_{\beta_0-\gamma_j}]) \\
				&=a-x_{\gamma_j}f_{D,\xi}([e_{\gamma_j},e_{\beta_0-\gamma_j}]) \\
                &= a-x_{\gamma_j}\cdot c_j \cdot \xi(\beta_0),
			\end{align*}
		where $c_j$ is the nonzero scalar such that $[e_{\gamma_j},e_{\beta_0-\gamma_j}]=c_je_{\beta_0}$, while
		\begin{equation*}
			a=\begin{cases}
				\xi(\beta),&\text{if }\beta_0-\gamma_j=\beta,\\
				0&\text{otherwise}.
			\end{cases}
		\end{equation*}
		Hence, all $x_{\gamma_j}$ are uniquely defined by $\mu$. Now, let $S_k$ be the symmetric group on $k$ letters. We obtain

			\begin{align*}
				(\exp x.f_{D,\xi})(e_\beta)&= f_{D,\xi}\Big(e_\beta+(-1)^k\cdot \frac{1}{k!}\underbrace{[x,[x,\ldots [x,e_\beta]\ldots ]]}_{\text{$k$ commutators}}\Big) \\
				&=\xi(\beta) + f_{D,\xi}\Big((-1)^k\cdot\frac{1}{k!}\Big[\sum\limits_{j=1}^kx_{\gamma_j}e_{\gamma_j},\Big[\sum\limits_{j=1}^kx_{\gamma_j}e_{\gamma_j},\ldots\! \Big[\sum\limits_{j=1}^kx_{\gamma_j}e_{\gamma_j}, e_\beta\Big]\!\ldots\!\Big]\Big]\Big) \\
				&=\xi(\beta) + (-1)^k\cdot\frac{1}{k!}\prod_{j=1}^kx_{\gamma_j}\cdot f_{D,\xi}\Big(\sum\limits_{\delta\in S_k}[e_{\gamma_{\delta(1)}},[e_{\gamma_{\delta(2)}},\ldots [e_{\gamma_{\delta(k)}},e_\beta]\ldots ]]\Big).
			\end{align*}
		
		Denote the second summand by $F$. Then $F$ is uniquely defined by $\mu$, because $x_{\gamma_j}$ and $\xi(\beta_0)$ are uniquely defined by $\mu$.
		If $\Omega_{D,\xi_1} = \Omega_{D,\xi_2}$ then there exist $x_1$, $x_2$, for which $\exp x_1.f_{D,\xi_1} = \exp x_2.f_{D,\xi_2}$. So, $(\exp x_1.f_{D,\xi_1})(e_\beta) = (\exp x_2.f_{D,\xi_2})(e_\beta)$, or equivalently, $\xi_1(\beta)+F = \xi_2(\beta)+F$, hence $\xi_1(\beta)= \xi_2(\beta)$, a contradiction.\end{proof}
	
		Now, we need a general construction, which can be applied to an arbitrary root system.
	Namely, let $\gt$, $\bt$ and $N=\exp(\nt)$ be as in the introduction. Let $\htt$ be the Cartan subalgebra of $\gt$ such that $\gt=\nt\oplus\htt\oplus\nt^-$, where $\nt^-$ is the nilradical of the Borel subalgebra opposite to $\bt$, and let $\Phi^-$ be the set of negative roots. Then the root vectors $e_{\alpha}$, $\alpha\in\Phi^-$, form a basis of $\nt^-$. Further, let $\alpha_1, \ldots, \alpha_n$ be the simple roots from $\Phi$, and $h_{\alpha_i}$, $1\leq i\leq n$, be a basis of $\htt$ such that $\{e_{\alpha},\alpha\in\Phi\}\cup\{h_{\alpha_i}, 1\leq i\leq n\}$ is a Chevalley basis of $\gt$.
	
	We fix a total order $\leq_t$ on this basis such that $e_{\alpha}<_th_{\alpha_i}<_te_{-\beta}$ for all $\alpha,\beta\in\Phi^+$, $1\leq i\leq n$, and $e_{\alpha}<_te_{\beta}$ if $\alpha,\beta\in\Phi$ and $\alpha>\beta$. This identifies $\glt(\gt)$ with the Lie algebra $\glt_{\dim\gt}(\Cp)$, and $\aad(\nt)$ with a subalgebra of the Lie algebra $\ut$ of all the strictly upper-triangular matrices of $\glt_{\dim\gt}(\Cp)$.
	
	Let $\GL(V)$ be the group of all invertible linear operators on a vector space $V$. Since we have fixed a basis for $\gt$, the group $\GL(\gt)$ can be identified with the group $\GL_{\dim\gt}(\Cp)$, and $\exp\aad(\nt)\cong N$ is identified with a subgroup of the group $U$ of all upper-triangular matrices from $\GL_{\dim\gt}(\Cp)$ with $1$'s on the diagonal. Furthermore, using the Killing form on $\gt$ and the trace form on $\glt(\gt)$, one can identify $\nt^*$ with the space $\nt_-=\langle e_{-\alpha}, \alpha\in\Phi^+\rangle_{\Cp}$ and $\ut^*$ with the space $\ut_-=\ut^T$, where the superscript $T$ denotes the transposed matrix. Under all these identifications, it is enough to check that the coadjoint $U$-orbits of the linear forms $\wt f_{D,\xi_1}$ and $\wt f_{D,\xi_2}$ are distinct. Here, given a map $\xi\colon D\to\Cp^{\times}$, we denote by $\wt f_{D,\xi}$ the matrix $$\wt f_{D,\xi}=\left(\aad\left(\sum\nolimits_{\beta\in D}\xi(\beta)e_{\beta}\right)\right)^T\in\ut_-\cong\ut^*.$$
	
	We will now study the matrix $f=\wt f_{D,\xi}$ in more detail. The rows and the columns of matrices from $\glt(\gt)$ are now indexed by the elements of the Chevalley basis fixed above. Given a matrix $x$ from $\glt(\gt)$ and two basis elements $a, b$, we will denote by $x_{a,b}$ the entry of $x$ lying in the $a$th row and the $b$th column. The following proposition was proved in \cite {IgnatyevShevchenko21}. For the reader's convenience, we reproduce the proof here, because our main technical tool used in the proof of Theorem~\ref{mtheo:F4} is based on similar ideas.

	\begin{proposition}[{\cite[Proposition~4.2]{IgnatyevShevchenko21}}]\label{prop:distinct_orbits} Let $\Phi$ be an irreducible root system\textup, and  $D$ be a non-singular rook placement in $\Phi^+$. Let $\beta_0$ be a root in $D$\textup, $\xi_1$ and $\xi_2$ be maps from $D$ to $\Cp^{\times}$ for which $\xi_1(\beta_0)\neq\xi_2(\beta_0)$. Assume that there exists a simple root $\alpha_0\in\Delta$ satisfying $(\alpha_0,\beta_0)\neq0$ and $(\alpha_0,\beta)=0$ for all $\beta\in D$ such that $\beta\nless\beta_0$. Then $\Omega_{D,\xi_1}\neq\Omega_{D,\xi_2}$.
 \end{proposition}
	\begin{proof}
Since $$\ad{e_{\beta_0}}(h_{\alpha_0})=[e_{\beta_0},h_{\alpha_0}]=-\dfrac{2(\alpha_0,\beta_0)}{(\alpha_0,\alpha_0)}e_{\beta_0},$$ we obtain $f_{h_{\alpha_0},e_{\beta_0}}=-\xi(\beta_0)\dfrac{2(\alpha_0,\beta_0)}{(\alpha_0,\alpha_0)}\neq0$. One may assume without loss of generality that $h_{\alpha_0}>_th_{\alpha_i}$ for all $\alpha_i\neq\alpha_0$. We claim that
		\begin{equation}\label{formula:condition_Andre}
			f_{h_{\alpha_0},e_{\alpha}}=f_{e_{-\gamma},e_{\beta_0}}=0\text{ for all }e_{\alpha}<_te_{\beta_0}\text{ and all }e_{-\gamma}, \alpha,\gamma\in\Phi^+.
		\end{equation}
		Indeed, if $\alpha\notin D$ then, evidently, $f_{h_{\alpha_0},e_{\alpha}}=0$. If $\alpha=\beta\in D$ and $e_{\beta}<_te_{\beta_0}$ then $\beta\nless\beta_0$, hence $$f_{h_{\alpha_0},e_{\beta}}=-\xi(\beta)\dfrac{2(\alpha_0,\beta)}{(\alpha_0,\alpha_0)}=0,$$ because $(\alpha_0,\beta)=0$. On the other hand, if $f_{e_{-\gamma},e_{\beta_0}}\neq0$ for some $\gamma\in\Phi^+$ then $\beta_0=\beta-\gamma$. This contradicts the condition $\beta_0\notin S(\beta)$.
		
		Thus, $(\wt f_{D,\xi_1})_{h_{\alpha_0},e_{\alpha}}$ and $(\wt f_{D,\xi_2})_{h_{\alpha_0},e_{\alpha}}$ are different nonzero scalars, and (\ref{formula:condition_Andre}) is satisfied both for $f=\wt f_{D,\xi_1}$ and for $f=\wt f_{D,\xi_2}$. Now it follows immediately from the proof of \cite[Proposition~3]{Andre95} (or from Remark~\ref{nota:Andre}) that the coadjoint $U$-orbits of these matrices are distinct, and, consequently, $\Omega_{D,\xi_1}\neq\Omega_{D,\xi_2}$, as required.\end{proof}
	
	Our main technical tool generalizes the proposition above in the following way. Fix an order $\{\beta_1, \ldots, \beta_m\}$ on $D$ and an order $\{\alpha_1, \ldots, \alpha_n\}$ on the simple roots in $\Phi^+$ such that $h_{\alpha_i}<_th_{\alpha_j}$ and $e_{\beta_i}<_te_{\beta_j}$ for $i<j$. Note that
	\begin{equation*}
		f_{h_{\alpha_i},e_{\beta_j}}=-\xi(\beta_j)\dfrac{2(\alpha_i,\beta_j)}{(\alpha_i,\alpha_i)}.
	\end{equation*}
	Given $J\in\{1, \ldots, m\}$ and $I\subset\{1, \ldots, n\}$ with $|I|=|J|$, denote by $\Delta_I^J(\xi)$ the minor of the matrix $f$ with the set of rows $\{h_{\alpha_i}, i\in I\}$ and the set of columns $\{e_{\beta_j}, j\in J\}$. Furthermore, let $\wt\Delta_I^J$ be the determinant of the matrix, which $(i,j)$-th element equals $p_{i,j}=\dfrac{2(\alpha_i,\beta_j)}{(\alpha_i,\alpha_i)}$, so that $\Delta_I^J(\xi)=\pm\prod_{j\in J}\xi(\beta_j)\wt\Delta_I^J$.
	
	
	\begin{proposition}\label{prop:distinct_orbits_gen}
 Assume that there exist an $m$-tuple $I=(i_1, \ldots, i_m)$ such that, for all $1\leq k\leq m$,  $\wt\Delta_{I_k}^{J_k}\neq0$, where $I_k=\{i_l \mid l\leq k, i_l\geq i_k\}$ and $J_k=\{j\mid i_j\in I_k\}$. Assume also that, for all $1\leq k\leq m$, $\wt\Delta_{I_l'}^{J_l'}=0$ for $l\notin \{i_1, \ldots, i_{k-1}\}$, $l>i_k$, where
 \[
    I_l'=\{l\}\cup\{i_s\mid s<k, i_s>l\}
    \quad \textup{ and } \quad
    J_l'=\{k\}\cup\{j\mid i_j\in I_l'\setminus\{l\}\}. 
\]
Let $\xi_1$ and $\xi_2$ be maps from $D$ to $\Cp^{\times}$. If $\xi_1\neq\xi_2$ then $\Omega_{D,\xi_1}\neq\Omega_{D, \xi_2}$.\end{proposition}
  \begin{proof}For simplicity, we denote $\Phi^+_1=\{\delta\mid e_{\delta}<_te_{\beta}\text{ for all }\beta\in D\}$ and $\Phi^+_2=\Phi^+\setminus(D\cup\Phi^+_1)$.
		
		First, note that $f_{e_{-\gamma},e_{\beta_j}}=0$ and $f_{h_{\alpha_i},e_{\delta}}=0$ for all $\gamma, \delta\in\Phi^+_1$, $\alpha_i\in\Delta$, $\beta_j\in D$. Indeed, $f_{e_{-\gamma},e_{\beta_j}}$ equals the coefficient of $e_{\beta_j}$ in the expression $\sum_{\beta\in D}\xi(\beta)[e_{\beta},e_{-\gamma}]$. But if this coefficient is nonzero then $\beta-\gamma=\beta_j$ for some $\beta\in D$, which contradicts the non-singularity of $D$. On the other hand, $f_{h_{\alpha_i},e_{\delta}}$ equals the coefficient of $e_{\delta}$ in the expression $\sum_{\beta\in D}\xi(\beta)[e_{\beta},h_{\alpha_i}]$ which is clearly zero, because $[e_{\beta},h_{\alpha_i}]$ is parallel to $e_{\beta}$ for each $\beta\in D$, while $\delta\notin D$. On the picture below we
		draw schematically the matrix $f$. Marks $\Phi^+_1$, $D$, $\Phi^+_2$, $\Delta$, $\Phi^-$ mean that the corresponding rows and columns of the matrix $f$ are indexed by $e_{\delta}$ for $\delta\in\Phi^+_1$, $e_{\beta_j}$ for $\beta_j\in D$, $e_{\gamma}$ for $\gamma\in\Phi^+_2$, $h_{\alpha_i}$ for $\alpha_i\in\Delta$, $e_{\alpha}$ for $\alpha\in\Phi^-$ respectively. We replaced by big zeroes the blocks $\Delta\times\Phi^+_1$ and $\Phi^-\times D$ filled in zero entries. The minors $\Delta_{I_l}^{J_l}$ are the determinants of submatrices of the grey block $\Delta\times D$.
		\begin{center}
			\includegraphics{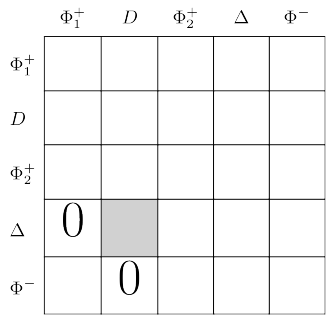}
		\end{center}
		
		The Lie algebra $\ut$ corresponds to the root system $A_{N-1}$, where $N=|\Phi|+\rk{\Phi}$. Let $\wt D_i$ be the subset of $A_{N-1}^+$ and $\wt\xi_i\colon\wt D_i\to\Cp^{\times}$ be the map such that $f_i=\wt f_{D,\xi_i}$ (as an element of $\ut^*$) belongs to the basic subvariety $\Ou_{\wt D_i,\wt\xi_i}$ of $\ut^*$ defined in Section~\ref{sect:Andre}, $i=1, 2$. Put $\Jo=\{e_{\alpha}, \alpha\in\Phi\}\cup\{h_{\alpha_i}, \alpha_i\in\Delta\}$. Each pair $(x,y)\in\Jo\times\Jo$ such that the $(x,y)$-th entry of $f_i$ lies under the diagonal corresponds to the unique root $\epsi_y-\epsi_x\in A_{N-1}^+$. We denote the inverse map from $A_{N-1}^+$ to $\Jo\times\Jo$ by $\tau$.
		
		Put $L=\{e_{\alpha}, \alpha\in\Phi^-\}\times\{e_{\delta, \delta\in\Phi^+_1}\}$ and $\wt D_i^L=\tau(\wt D_i)\cap L$. According to Andr\'e's theory, we may assume without loss of generality that $\wt D_1^L=\wt D_2^L$ and $\wt\xi_1(\epsi_y-\epsi_x)=\wt\xi_2(\epsi_y-\epsi_x)$ for each $(x,y)\in\wt D_i^L$ (if not, then $\Ou_{\wt D_1,\wt\xi_1}\neq\Ou_{\wt D_2,\wt\xi_2}$ and, consequently, $\Omega_{D,\xi_1}\neq\Omega_{D,\xi_2}$.) We will prove that $\xi_1(\beta_j)=\xi_2(\beta_j)$ for all $1\leq j\leq m$ by induction on $j$. The case $j=0$ (with $I_0=\varnothing$) can be considered as an evident inductive base case.
		
		Let $j\geq1$. Note that each $\beta_l$, $1\leq l\leq m$, belongs to $\wt D_i$, and the intersection of $\tau(\wt D_i)$ with $\{h_{\alpha_i}, i\in\Delta\}\times\{e_{\beta}, \beta\in D\}$ (i.e., with the ``grey'' area) coincides with $\{(h_{\alpha_{i_l}},\beta_l)\}_{l=1}^m$ (this follows immediately from Remark~\ref{nota:Andre}). Furthermore, recall the notion of $\Delta_{\wt\alpha}^{\wt D_i}(f_i)$ for $\wt\alpha\in A_{N-1}^+$ from Section~\ref{sect:Andre}, where $\wt f_i$ is considered as an element of $\ut^*$. It also follows from Remark~\ref{nota:Andre} that, for each $l$ from $1$ to $m$ and for $i=1, 2$,
		\begin{equation*}
			\Delta_{\tau^{-1}(h_{\alpha_{i_l}},e_{\beta_l})}^{\wt D_i}=\pm\prod\nolimits_{\wt\alpha\in\wt D_i, \tau(\wt\alpha)\in\wt D_i^L}\wt\xi_i(\wt\alpha)\Delta_{I_l}^{J_l}=
			\mathrm{const}_{i,l}\wt\Delta_{I_l}^{J_l},
		\end{equation*}
		where $\mathrm{const}_{i,l}$ is a scalar depending only on $D_{i,l}=\wt D_i^L\cup\{(h_{\alpha_{i_s}}, e_{\beta_s}), s<l, h_{\alpha_{i_s}}>_th_{\alpha_{i_l}}\}$ and on $\restr{\wt\xi_i}{D_{i,l}}$. By the inductive assumption, $D_{1,j}=D_{2,j}$ and $\restr{\wt\xi_2}{D_{i,j}}=\restr{\wt\xi_2}{D_{2,j}}$. We conclude that $\wt\xi_1(\beta_j)=\wt\xi_2(\beta_j)$, as required and the proof is complete.\end{proof}
	
\begin{remark}It follows from the conditions of Proposition~\ref{prop:distinct_orbits_gen} that if such an $m$-tuple $I$ exists then it is unique.
 \end{remark}
	
	From now on, let $\Phi=F_4$. Recall that the set $\Delta$ of simple roots can be identified with the following subset of $\Rp^4$:
	\begin{equation*}
		\Delta=\{\alpha_1=\epsi_2-\epsi_3, \alpha_2=\epsi_3-\epsi_4, \alpha_3=\epsi_4, \alpha_4=\dfrac{1}{2}(\epsi_1-\epsi_2-\epsi_3-\epsi_4)\}.
	\end{equation*}
	Here $\{\epsi_i\}_{i=1}^4$ is the standard basis of $\Rp^4$ (with the standard inner product). The set of positive roots is as follows:
	
		\begin{align*}
			\Phi^+\! =\{&\alpha_1, ~\alpha_2, ~\alpha_3, ~\alpha_4, ~\alpha_1+\alpha_2, ~\alpha_2+\alpha_3, ~\alpha_2+2\alpha_3, ~\alpha_3+\alpha_4, ~\alpha_1+\alpha_2+\alpha_3,\\
            &\alpha_1+\alpha_2+2\alpha_3,~\alpha_1+2\alpha_2+2\alpha_3, ~\alpha_2+\alpha_3+\alpha_4, ~\alpha_2+2\alpha_3+\alpha_4, ~\alpha_2+2\alpha_3+2\alpha_4,\\
            &\alpha_1+\alpha_2+\alpha_3+\alpha_4,~\alpha_1+\alpha_2+2\alpha_3+\alpha_4, ~\alpha_1+\alpha_2+2\alpha_3+2\alpha_4,\\
            &\alpha_1+2\alpha_2+2\alpha_3+\alpha_4,~\alpha_1+2\alpha_2+2\alpha_3+2\alpha_4 \alpha_1+2\alpha_2+3\alpha_3+\alpha_4,\\
            &\alpha_1+2\alpha_2+3\alpha_3+2\alpha_4,~\alpha_1+2\alpha_2+4\alpha_3+2\alpha_4,~\alpha_1+3\alpha_2+4\alpha_3+2\alpha_4,\\
            &2\alpha_1+3\alpha_2+4\alpha_3+2\alpha_4\}.
		\end{align*}

	We will apply Proposition~\ref{prop:distinct_orbits_gen} above to the following rook placements.
	
	\begin{proposition}\label{prop:prop_apply}
 Let $\Phi=F_4$\textup, and $D=\{\beta_1, \ldots, \beta_m\}$ be one of rook placements from the table below. Then the orders on $\Delta$ and $D$ and the sets $I_j=\{i_1, \ldots, i_j\}$\textup, $1\leq j\leq m$\textup, from the table below satisfy the conditions of Proposition~\textup{\ref{prop:distinct_orbits_gen}}.
		{\small
  \begin{longtable}{||l|l|l|l||}
			\hline
			&Rook placement $D$&Order on $\Delta$&$i_1, \ldots, i_m$\\
			\hline\hline
			\textup{1}&
			\hspace{-0.3cm}\begin{tabular}[t]{l}
				$\beta_1=\alpha_1+3\alpha_2+4\alpha_3+2\alpha_4$, $\beta_2=\alpha_1+\alpha_2+2\alpha_3+2\alpha_4$,\\
				$\beta_3=\alpha_1+\alpha_2+2\alpha_3$
			\end{tabular}&
			$\alpha_1, \alpha_4, \alpha_2, \alpha_3$&
			$3, 2, 1$
			\\
			\hline
			\textup{2}&
			\hspace{-0.3cm}\begin{tabular}[t]{l}
				$\beta_1=\alpha_1+\alpha_2+2\alpha_3+2\alpha_4$, $\beta_2=\alpha_1+\alpha_2+\alpha_3$
			\end{tabular}&
			$\alpha_2, \alpha_1, \alpha_4, \alpha_3$&
			$3, 2$
			\\
			\hline
			\textup{3}&
			\hspace{-0.3cm}\begin{tabular}[t]{l}
				$\beta_1=\alpha_1+3\alpha_2+4\alpha_3+2\alpha_4$, $\beta_2=\alpha_1+\alpha_2+2\alpha_3+2\alpha_4$,\\
				$\beta_3=\alpha_1+\alpha_2+2\alpha_3$,
				$\beta_4=\alpha_1+\alpha_2$
			\end{tabular}&
			$\alpha_1, \alpha_2, \alpha_4, \alpha_3$&
			$2, 3, 4, 1$
			\\
			\hline
			\textup{4}&
			\hspace{-0.3cm}\begin{tabular}[t]{l}
				$\beta_1=\alpha_1+\alpha_2+2\alpha_3+2\alpha_4$, $\beta_2=\alpha_1+\alpha_2+2\alpha_3$,\\
				$\beta_3=\alpha_1+\alpha_2$
			\end{tabular}&
			$\alpha_1, \alpha_2, \alpha_4, \alpha_3$&
			$3, 4, 2$
			\\
			\hline
			\textup{5}&
			\hspace{-0.3cm}\begin{tabular}[t]{l}
				$\beta_1=\alpha_1+3\alpha_2+4\alpha_3+2\alpha_4$, $\beta_2=\alpha_1+\alpha_2+2\alpha_3$,\\
				$\beta_3=\alpha_1+\alpha_2$
			\end{tabular}&
			$\alpha_1, \alpha_2, \alpha_3, \alpha_4$&
			$2, 4, 3$
			\\
			\hline
			\textup{6}&
			\hspace{-0.3cm}\begin{tabular}[t]{l}
				$\beta_1=\alpha_1+3\alpha_2+4\alpha_3+2\alpha_4$, $\beta_2=\alpha_1+\alpha_2+2\alpha_3+\alpha_4$,\\
				$\beta_3=\alpha_1+\alpha_2$
			\end{tabular}&
			$\alpha_1, \alpha_2, \alpha_3, \alpha_4$&
			$2, 3, 1$
			\\
			\hline
			\textup{7}&
			\hspace{-0.3cm}\begin{tabular}[t]{l}
				$\beta_1=\alpha_1+\alpha_2+2\alpha_3+\alpha_4$, $\beta_2=\alpha_1+\alpha_2$
			\end{tabular}&
			$\alpha_1, \alpha_2, \alpha_3, \alpha_4$&
			$3, 2$
			\\
			\hline
			\textup{8}&
			\hspace{-0.3cm}\begin{tabular}[t]{l}
				$\beta_1=\alpha_1+3\alpha_2+4\alpha_3+2\alpha_4$, $\beta_2=\alpha_1+2\alpha_2+2\alpha_3+2\alpha_4$,\\
				$\beta_3=\alpha_1+2\alpha_2+2\alpha_3$,
				$\beta_4=\alpha_1$
			\end{tabular}&
			$\alpha_2, \alpha_3, \alpha_4, \alpha_1$&
			$2, 3, 1, 4$
			\\
			\hline
			\textup{9}&
			\hspace{-0.3cm}\begin{tabular}[t]{l}
				$\beta_1=\alpha_1+2\alpha_2+4\alpha_3+2\alpha_4$, $\beta_2=\alpha_1+2\alpha_2+2\alpha_3+\alpha_4$,\\
				$\beta_3=\alpha_1$
			\end{tabular}&
			$\alpha_3, \alpha_2, \alpha_1, \alpha_4$&
			$2, 1, 3$
			\\
			\hline
			\textup{10}&
			\hspace{-0.3cm}\begin{tabular}[t]{l}
				$\beta_1=\alpha_1+2\alpha_2+2\alpha_3+2\alpha_4$, $\beta_2=\alpha_3+\alpha_4$
			\end{tabular}&
			$\alpha_4, \alpha_3, \alpha_2, \alpha_1$&
			$3, 2$
			\\
			\hline
			\textup{11}&
			\hspace{-0.3cm}\begin{tabular}[t]{l}
				$\beta_1=2\alpha_1+3\alpha_2+4\alpha_3+2\alpha_4$, $\beta_2=\alpha_2+\alpha_3+\alpha_4$,\\
				$\beta_3=\alpha_2+2\alpha_3$
			\end{tabular}&
			$\alpha_1, \alpha_3, \alpha_4, \alpha_2$&
			$1, 4, 3$
			\\
			\hline
			\textup{12}&
			\hspace{-0.3cm}\begin{tabular}[t]{l}
				$\beta_1=\alpha_1+2\alpha_2+2\alpha_3+\alpha_4$, $\beta_2=\alpha_2+2\alpha_3+2\alpha_4$,\\
				$\beta_3=\alpha_2+2\alpha_3$
			\end{tabular}&
			$\alpha_1, \alpha_4, \alpha_3, \alpha_2$&
			$4, 2, 3$
			\\
			\hline
			\textup{13}&
			\hspace{-0.3cm}\begin{tabular}[t]{l}
				$\beta_1=\alpha_1+3\alpha_2+4\alpha_3+2\alpha_4$, $\beta_2=\alpha_1+\alpha_2+\alpha_3+\alpha_4$,\\
				$\beta_3=\alpha_1+\alpha_2+2\alpha_3$
			\end{tabular}&
			$\alpha_3, \alpha_4, \alpha_1, \alpha_2$&
			$4, 3, 2$
			\\
			\hline
			\textup{14}&
			\hspace{-0.3cm}\begin{tabular}[t]{l}
				$\beta_1=\alpha_1+3\alpha_2+4\alpha_3+2\alpha_4$, $\beta_2=\alpha_1+\alpha_2+2\alpha_3$,\\
				$\beta_3=\alpha_3+\alpha_4$
			\end{tabular}&
			$\alpha_1, \alpha_2, \alpha_3, \alpha_4$&
			$2, 4, 3$
			\\
			\hline
			\textup{15}&
			\hspace{-0.3cm}\begin{tabular}[t]{l}
				$\beta_1=\alpha_1+3\alpha_2+4\alpha_3+2\alpha_4$, $\beta_2=\alpha_1+\alpha_2$, $\beta_3=\alpha_4$
			\end{tabular}&
			$\alpha_1, \alpha_2, \alpha_3, \alpha_4$&
			$2, 3, 4$
			\\
			\hline
			\textup{16}&
			\hspace{-0.3cm}\begin{tabular}[t]{l}
				$\beta_1=\alpha_1+2\alpha_2+3\alpha_3+\alpha_4$, $\beta_2=\alpha_1+\alpha_2+2\alpha_3+2\alpha_4$,\\
				$\beta_3=\alpha_1+\alpha_2$
			\end{tabular}&
			$\alpha_2, \alpha_3, \alpha_1, \alpha_2$&
			$2, 4, 3$
			\\
			\hline
			\textup{17}&
			\hspace{-0.3cm}\begin{tabular}[t]{l}
				$\beta_1=\alpha_1+\alpha_2+2\alpha_3+2\alpha_4$, $\beta_2=\alpha_2+\alpha_3+\alpha_4$,\\
				$\beta_3=\alpha_1+\alpha_2$
			\end{tabular}&
			$\alpha_4, \alpha_1, \alpha_2, \alpha_3$&
			$3, 4, 2$
			\\
			\hline
			\textup{18}&
			\hspace{-0.3cm}\begin{tabular}[t]{l}
				$\beta_1=\alpha_1+2\alpha_2+4\alpha_3+2\alpha_4$, $\beta_2=\alpha_1+\alpha_2+\alpha_3+\alpha_4$,\\
				$\beta_3=\alpha_1+2\alpha_2+2\alpha_3$
			\end{tabular}&
			$\alpha_3, \alpha_2, \alpha_4, \alpha_1$&
			$2, 4, 3$
			\\
			\hline
			\textup{19}&
			\hspace{-0.3cm}\begin{tabular}[t]{l}
				$\beta_1=\alpha_1+2\alpha_2+4\alpha_3+2\alpha_4$, $\beta_2=\alpha_2+2\alpha_3+\alpha_4$,\\
				$\beta_3=\alpha_1+2\alpha_2+2\alpha_3$
			\end{tabular}&
			$\alpha_1, \alpha_2, \alpha_4, \alpha_3$&
			$2, 4, 3$
			\\
			\hline
			\textup{20}&
			\hspace{-0.3cm}\begin{tabular}[t]{l}
				$\beta_1=\alpha_1+2\alpha_2+4\alpha_3+2\alpha_4$, $\beta_2=\alpha_1+2\alpha_2+2\alpha_3+2\alpha_4$
			\end{tabular}&
			$\alpha_2, \alpha_3, \alpha_1, \alpha_4$&
			$2, 1$
			\\
			\hline
			\textup{21}&
			\hspace{-0.3cm}\begin{tabular}[t]{l}
				$\beta_1=\alpha_1+2\alpha_2+2\alpha_3+2\alpha_4$, $\beta_2=\alpha_3+\alpha_4$, $\beta_3=\alpha_1$
			\end{tabular}&
			$\alpha_3, \alpha_4, \alpha_2, \alpha_1$&
			$3, 2, 4$
			\\
			\hline
			\textup{22}&
			\hspace{-0.3cm}\begin{tabular}[t]{l}
				$\beta_1=\alpha_1+\alpha_2+2\alpha_3+\alpha_4$, $\beta_2=\alpha_2+2\alpha_3+2\alpha_4$,\\
				$\beta_3=\alpha_2+2\alpha_3$
			\end{tabular}&
			$\alpha_2, \alpha_1, \alpha_3, \alpha_4$&
			$3, 4, 2$
			\\
			\hline
			\textup{23}&
			\hspace{-0.3cm}\begin{tabular}[t]{l}
				$\beta_1=\alpha_1+2\alpha_2+2\alpha_3+\alpha_4$, $\beta_2=\alpha_1+\alpha_2+2\alpha_3+2\alpha_4$,\\
				$\beta_3=\alpha_1+\alpha_2+2\alpha_3$
			\end{tabular}&
			$\alpha_2, \alpha_1, \alpha_3, \alpha_4$&
			$3, 4, 2$
			\\
			\hline
			\textup{24}&
			\hspace{-0.3cm}\begin{tabular}[t]{l}
				$\beta_1=\alpha_1+2\alpha_2+3\alpha_3+2\alpha_4$, $\beta_2=\alpha_1+\alpha_2+2\alpha_3$,\\
				$\beta_3=\alpha_1+\alpha_2$ 
			\end{tabular}&
			$\alpha_2, \alpha_1, \alpha_3, \alpha_4$&
			$4, 3, 2$
			\\
			\hline
		\end{longtable}}
	\end{proposition}
 \begin{proof}The proof is case-by-case and is completely straightforward. As an example, consider the 17th rook placement $D$.
		
		Clearly, the root $\beta_3$ (respectively, $\beta_1$) is orthogonal to the unique simple root, namely, to $\alpha_3$ (respectively, to $\alpha_4$). There are no simple roots orthogonal to $\beta_2$. Write out the minor of the matrix $f$, which rows correspond to  $h_{\alpha_i}$, $\alpha_i\in\Delta$, and columns correspond to $e_{\beta_j}$, $\beta_j\in D$. Recall the notion of $p_{i,j}$ introduced before Proposition~\ref{prop:distinct_orbits_gen}.
		\begin{equation*}
			\begin{matrix}
				& \beta_3& \beta_2& \beta_1\\
				\alpha_4& p_{1,1}& p_{1,2}& 0\\
				\alpha_1& p_{2,1}& p_{2,2}& p_{2,3}\\
				\alpha_2& p_{3,1}& p_{3,2}& p_{3,3}\\
				\alpha_3& 0& p_{4,2}& p_{4,3}
			\end{matrix}
		\end{equation*}
		
		Obviously, $\wt\Delta_{\{4\}}^{\{1\}}=|p_{4,1}|=0$, while, for $i_1=3$, $$\wt\Delta_{I_1}^{J_1}=|p_{3,1}|=\dfrac{2(\alpha_2,\beta_3)}{(\alpha_2,\alpha_2)}=-1\neq0.$$ Hence, in fact we have the only possibility for $i_1$: $i_1=3$. Next, for $i_2=4$, one has $$\wt\Delta_{I_2}^{J_2}=|p_{4,2}|=\dfrac{2(\alpha_3,\beta_2)}{(\alpha_3,\alpha_3)}=-1\neq0.$$ Therefore, we have to put $i_2=4$. Finally, for $i_3=2$, we obtain
			\begin{align*}
				\wt\Delta_{I_3}^{J_3}=
				\begin{vmatrix}
					p_{2,1}& p_{2,2}& p_{2,3}\\
					p_{3,1}& p_{3,2}& p_{3,3}\\
					0& p_{4,2}& p_{4,3}
				\end{vmatrix}
				&=\dfrac{8}{(\alpha_1,\alpha_1)(\alpha_2,\alpha_2)(\alpha_3,\alpha_3)}\begin{vmatrix}
					(\alpha_1,\beta_3)& (\alpha_1,\beta_2)& (\alpha_1,\beta_1)\\
					(\alpha_2,\beta_3)& (\alpha_2,\beta_2)& (\alpha_2,\beta_1)\\
					0& (\alpha_3,\beta_2)& (\alpha_3,\beta_1)
				\end{vmatrix}\\
				&=2
				\begin{vmatrix}
					1& -1& 1\\
					-1& 1& 1\\
					0& -\frac{1}{2}& -1
				\end{vmatrix}
				=4\neq0.
			\end{align*}
		Thus, there is the only candidate for $i_3$: $i_3=2$. It is easy to check that the sequence $(3, 4, 2)$ satisfies the conditions of Proposition~\ref{prop:distinct_orbits_gen}.
		
		All other rook placements from the table above can be considered similarly.\end{proof}
 
	We are now ready to prove our second main result, Theorem~\ref{mtheo:F4}, which claims that, for a non-singular orthogonal rook placement $D\subset F_4^+$ and two distinct maps $\xi_1$, $\xi_2$ from $D$ to $\Cp^{\times}$, the associated coadjoint orbits $\Omega_{D,\xi_1}$, $\Omega_{D,\xi_2}$ do not coincide.
	
	\begin{proof}[Proof of Theorem~\ref{mtheo:F4}] This proof is based on a case-by-case analysis. Namely, we split the rook placements in $F_4^+$ into several ``classes'' and then apply Lemma~\ref{lemm:F4_submax} and Propositions~\ref{prop:distinct_orbits}, \ref{prop:distinct_orbits_gen} to these classes. We start with the maximal (possibly, singular) rook placements. It is easy to check that there are 24 maximal rook placements in $F_4^+$:
 \begingroup
 \allowdisplaybreaks
	\begin{align*}
			D_1=\{&\beta_1=\alpha_1+2\alpha_2+3\alpha_3+2\alpha_4, \beta_2=\alpha_1+\alpha_2+\alpha_3, \beta_3=\alpha_2+\alpha_3, \beta_4=\alpha_3\},\\
			D_2=\{&\beta_1=\alpha_1+2\alpha_2+3\alpha_3+2\alpha_4, \beta_2=\alpha_1+\alpha_2+\alpha_3, \beta_3=\alpha_2+2\alpha_3, \beta_4=\alpha_2\},\\
			D_3=\{&\beta_1=\alpha_1+2\alpha_2+3\alpha_3+2\alpha_4, \beta_2=\alpha_1+2\alpha_2+2\alpha_3, \beta_3=\alpha_3, \beta_4=\alpha_1\},\\
			D_4=\{&\beta_1=\alpha_1+2\alpha_2+4\alpha_3+2\alpha_4, \beta_2=\alpha_1+2\alpha_2+2\alpha_3+2\alpha_4, \beta_3=\alpha_1+\alpha_2+\alpha_3,\\
                &\beta_4=\alpha_2+\alpha_3\},\\
			D_5=\{&\beta_1=2\alpha_1+3\alpha_2+4\alpha_3+2\alpha_4, \beta_2=\alpha_2+2\alpha_3+2\alpha_4, \beta_3=\alpha_2+\alpha_3, \beta_4=\alpha_3\},\\
			D_6=\{&\beta_1=2\alpha_1+3\alpha_2+4\alpha_3+2\alpha_4, \beta_2=\alpha_2+2\alpha_3+2\alpha_4, \beta_3=\alpha_2+2\alpha_3, \beta_4=\alpha_2\},\\
			D_7=\{&\beta_1=2\alpha_1+3\alpha_2+4\alpha_3+2\alpha_4, \beta_2=\alpha_2+2\alpha_3+\alpha_4, \beta_3=\alpha_4, \beta_4=\alpha_2\},\\
			D_8=\{&\beta_1=\alpha_1+2\alpha_2+3\alpha_3+\alpha_4, \beta_2=\alpha_1+\alpha_2+\alpha_3+\alpha_4, \beta_3=\alpha_2+2\alpha_3+2\alpha_4,\\     
                &\beta_4=\alpha_2\},\\
			D_9=\{&\beta_1=\alpha_1+2\alpha_2+3\alpha_3+\alpha_4, \beta_2=\alpha_1+\alpha_2+\alpha_3+\alpha_4, \beta_3=\alpha_2+\alpha_3+\alpha_4,\\
                &\beta_4=\alpha_3+\alpha_4\},\\
			D_{10}=\{&\beta_1=\alpha_1+2\alpha_2+2\alpha_3+\alpha_4, \beta_2=\alpha_1+\alpha_2+2\alpha_3+\alpha_4, \beta_3=\alpha_2+2\alpha_3+\alpha_4,\\
                &\beta_4=\alpha_4\},\\
			D_{11}=\{&\beta_1=\alpha_1+2\alpha_2+3\alpha_3+2\alpha_4, \beta_2=\alpha_1+\alpha_2+2\alpha_3, \beta_3=\alpha_2+\alpha_3, \beta_4=\alpha_1+\alpha_2\},\\
			D_{12}=\{&\beta_1=\alpha_1+3\alpha_2+4\alpha_3+2\alpha_4, \beta_2=\alpha_1+\alpha_2+2\alpha_3+2\alpha_4, \beta_3=\alpha_1+\alpha_2+\alpha_3, \\
                &\beta_4=\alpha_3\},\\
			D_{13}=\{&\beta_1=2\alpha_1+3\alpha_2+4\alpha_3+2\alpha_4, \beta_2=\alpha_2+\alpha_3+\alpha_4, \beta_3=\alpha_3+\alpha_4, \beta_4=\alpha_2+2\alpha_3\},\\
			D_{14}=\{&\beta_1=\alpha_1+2\alpha_2+2\alpha_3+\alpha_4, \beta_2=\alpha_1+\alpha_2+2\alpha_3+\alpha_4, \beta_3=\alpha_2+2\alpha_3+2\alpha_4,\\
                &\beta_4=\alpha_2+2\alpha_3\},\\
			D_{15}=\{&\beta_1=\alpha_1+3\alpha_2+4\alpha_3+2\alpha_4, \beta_2=\alpha_1+\alpha_2+2\alpha_3+2\alpha_4, \beta_3=\alpha_1+\alpha_2+2\alpha_3,\\
                &\beta_4=\alpha_1+\alpha_2\},\\
			D_{16}=\{&\beta_1=\alpha_1+3\alpha_2+4\alpha_3+2\alpha_4, \beta_2=\alpha_1+\alpha_2+\alpha_3+\alpha_4, \beta_3=\alpha_1+\alpha_2+2\alpha_3,\\
                &\beta_4=\alpha_3+\alpha_4\},\\
			D_{17}=\{&\beta_1=\alpha_1+3\alpha_2+4\alpha_3+2\alpha_4,\beta_2=\alpha_1+\alpha_2+2\alpha_3+\alpha_4, \beta_3=\alpha_1+\alpha_2, \beta_4=\alpha_4\},\\
			D_{18}=\{&\beta_1=\alpha_1+2\alpha_2+2\alpha_3+\alpha_4, \beta_2=\alpha_1+\alpha_2+2\alpha_3+2\alpha_4, \beta_3=\alpha_2+2\alpha_3+\alpha_4,\\
			&\beta_4=\alpha_1+\alpha_2+2\alpha_3\},\\
			D_{19}=\{&\beta_1=\alpha_1+2\alpha_2+3\alpha_3+\alpha_4, \beta_2=\alpha_1+\alpha_2+2\alpha_3+2\alpha_4, \beta_3=\alpha_2+\alpha_3+\alpha_4,\\
                &\beta_4=\alpha_1+\alpha_2\},\\
			D_{20}=\{&\beta_1=\alpha_1+2\alpha_2+4\alpha_3+2\alpha_4, \beta_2=\alpha_1+2\alpha_2+2\alpha_3+2\alpha_4, \beta_3=\alpha_1+2\alpha_2+2\alpha_3,\\
                &\beta_4=\alpha_1\},\\
                & \\
			D_{21}=\{&\beta_1=\alpha_1+2\alpha_2+4\alpha_3+2\alpha_4, \beta_2=\alpha_1+\alpha_2+\alpha_3+\alpha_4, \beta_3=\alpha_1+2\alpha_2+2\alpha_3,\\
			&\beta_4=\alpha_2+\alpha_3+\alpha_4\},\\
			D_{22}=\{&\beta_1=\alpha_1+2\alpha_2+4\alpha_3+2\alpha_4, \beta_2=\alpha_1+2\alpha_2+2\alpha_3+\alpha_4, \beta_3=\alpha_4, \beta_4=\alpha_1\},\\
			D_{23}=\{&\beta_1=\alpha_1+2\alpha_2+2\alpha_3+2\alpha_4, \beta_2=\alpha_1+\alpha_2+2\alpha_3+\alpha_4, \beta_3=\alpha_2+2\alpha_3+\alpha_4,\\
			&\beta_4=\alpha_1+2\alpha_2+2\alpha_3\},\\
			D_{24}=\{&\beta_1=\alpha_1+2\alpha_2+3\alpha_3+\alpha_4, \beta_2=\alpha_1+2\alpha_2+2\alpha_3+2\alpha_4, \beta_3=\alpha_3+\alpha_4, \beta_4=\alpha_1\}.
	\end{align*}

    \endgroup
	The first root $\beta_1$ is maximal among all roots in each of these rook placements. In the rook placements $$D_8, D_{14}, D_{18}, D_{19}, D_{24}$$ the second root $\beta_2$ is maximal, too. As we mentioned above, if $D$ is a subset of $D_i$ containing a maximal root $\beta$ from $D_i$ and $\Omega_{D,\xi_1}=\Omega_{D,\xi_2}$ then $\xi_1(\beta)=\xi_2(\beta)$ (here $\xi_1$, $\xi_2$ are maps from $D$ to $\Cp^{\times}$).
	
	Next, it is straightforward to check that the following maximal rook placements $D_i$ (together with a simple root $\alpha_0$ and a distinguished root $\beta_0\in D_i$) satisfy the conditions of Proposition~\ref{prop:distinct_orbits}, except the non-singularity of $D_i$:
	
	\begin{longtable}{lll}
		$D_2, \beta_0=\beta_3, \alpha_0=\alpha_1$;&
		$D_2, \beta_0=\beta_4, \alpha_0=\alpha_2$;&
		$D_3, \beta_0=\beta_2, \alpha_0=\alpha_2$;\\
		$D_3, \beta_0=\beta_4, \alpha_0=\alpha_1$;&
		$D_4, \beta_0=\beta_2, \alpha_0=\alpha_4$;&
		$D_4, \beta_0=\beta_3, \alpha_0=\alpha_1$;\\
		$D_5, \beta_0=\beta_2, \alpha_0=\alpha_4$;&
		$D_5, \beta_0=\beta_3, \alpha_0=\alpha_2$;&
		$D_6, \beta_0=\beta_2, \alpha_0=\alpha_4$;\\
		$D_6, \beta_0=\beta_3, \alpha_0=\alpha_3$;&
		$D_6, \beta_0=\beta_4, \alpha_0=\alpha_2$;&
		$D_7, \beta_0=\beta_2, \alpha_0=\alpha_3$;\\
		$D_7, \beta_0=\beta_4, \alpha_0=\alpha_2$;&
		$D_8, \beta_0=\beta_4, \alpha_0=\alpha_2$;&
		$D_{11}, \beta_0=\beta_2, \alpha_0=\alpha_1$;\\
		$D_{12}, \beta_0=\beta_2, \alpha_0=\alpha_2$;&
		$D_{13}, \beta_0=\beta_2, \alpha_0=\alpha_2$;&
		$D_{15}, \beta_0=\beta_2, \alpha_0=\alpha_4$;\\
		$D_{15}, \beta_0=\beta_3, \alpha_0=\alpha_3$;&
		$D_{17}, \beta_0=\beta_2, \alpha_0=\alpha_3$;&
		$D_{20}, \beta_0=\beta_2, \alpha_0=\alpha_4$;\\
		$D_{20}, \beta_0=\beta_4, \alpha_0=\alpha_1$;&
		$D_{21}, \beta_0=\beta_2, \alpha_0=\alpha_1$;&
		$D_{22}, \beta_0=\beta_4, \alpha_0=\alpha_1$;\\
		$D_{23}, \beta_0=\beta_2, \alpha_0=\alpha_1$;&
		$D_{24}, \beta_0=\beta_4, \alpha_0=\alpha_1$.&\\
	\end{longtable}
	This implies that if $D$ is a non-singular rook placement contained in one of these maximal rook placements and containing the root $\beta_0$, then $D$, $\beta_0$, $\alpha_0$ satisfy the conditions of Proposition~\ref{prop:distinct_orbits}. Hence, if $\Omega_{D,\xi_1}=\Omega_{D,\xi_2}$ then $\xi_1(\beta_0)=\xi_2(\beta_0)$.
	
	Now, let $D$ be a non-singular subset of one of the rook placements $D_1, \ldots, D_{10}$. Assume that $D\subset D_1$. Note that $\beta_i\in S(\beta_j)$ for all $2\leq i\leq4$ and $1\leq j<i$. Hence, $|D|=1$, and there is nothing to prove. If $D\subset D_2$ contains $\beta_1$ then $\beta_2\notin D$, because $\beta_2\in S(\beta_1)$. On the other hand, if $\beta_1\notin D$ and $\beta_2\in D$ then $\beta_2$ is maximal in $D$. For $\beta_3$, $\beta_4$ see the previous paragraph. Another example: assume that $D\subset D_3$. If $\beta_1\in D$ then $\beta_3\notin D$, because $\beta_3\in S(\beta_1)$. If $\beta_1\notin D$ and $\beta_3\in D$ then $D$, $\beta_0=\beta_3$, $\alpha_0=\alpha_3$ satisfy the condition of Proposition~\ref{prop:distinct_orbits}. For the roots $\beta_2$, $\beta_4$, see the previous paragraph. All other rook placements $D_4, \ldots, D_{10}$ can be considered in a similar way.
	
	Most of the remaining rook placements (i.e., non-singular subsets of $D_{11}, \ldots, D_{24}$) can be considered by completely similar arguments. The exceptions are the 24 rook placements from Proposition~\ref{prop:prop_apply} and the 8 following rook placements:
		\begin{align*}
			&D_{25}=\{\beta
			_1=\alpha_1+\alpha_2+2\alpha_3, \beta_2=\alpha_2+\alpha_3, \beta_3=\alpha_1+\alpha_2\};\\
			&D_{26}=\{\beta
			_1=\alpha_1+\alpha_2+2\alpha_3, \beta_2=\alpha_2+\alpha_3\};\\
			&D_{27}=\{\beta
			_1=\alpha_1+\alpha_2+2\alpha_3, \beta_2=\alpha_1+\alpha_2\};\\ &D_{28}=\{\beta_1=\alpha_1+\alpha_2+2\alpha_3+2\alpha_4, \beta_2=\alpha_2+2\alpha_3+\alpha_4, \beta_3=\alpha_1+\alpha_2+2\alpha_3\};\\ &D_{29}=\{\beta_1=\alpha_1+2\alpha_2+2\alpha_3+2\alpha_4, \beta_2=\alpha_1+2\alpha_2+2\alpha_3, \beta_3=\alpha_1\};\\ &D_{30}=\{\beta_1=\alpha_1+2\alpha_2+2\alpha_3+2\alpha_4, \beta_2=\alpha_1+2\alpha_2+2\alpha_3\};\\ &D_{31}=\{\beta_1=\alpha_1+2\alpha_2+2\alpha_3+2\alpha_4, \beta_2=\alpha_1+\alpha_2+2\alpha_3+\alpha_4, \beta_3=\alpha_1+2\alpha_2+2\alpha_3\};\\ &D_{32}=\{\beta_1=\alpha_1+2\alpha_2+2\alpha_3+2\alpha_4, \beta_2=\alpha_2+2\alpha_3+\alpha_4, \beta_3=\alpha_1+2\alpha_2+2\alpha_3\}.
		\end{align*}
	Proposition~\ref{prop:distinct_orbits_gen} completes the proof for 24 rook placements from Proposition~\ref{prop:prop_apply}. For the rook placements $D_{25}, \ldots, D_{32}$, one can apply Lemma~\ref{lemm:F4_submax} with $\beta=\beta_2$ or $\beta_3$ for $D_{25}$, $\beta=\beta_2$ for $D_{26}$, $\beta=\beta_2$ for $D_{27}$, $\beta=\beta_2$ or $\beta_3$ for $D_{28}$, $\beta=\beta_2$ for $D_{29}$, $\beta=\beta_2$ for $D_{30}$, $\beta=\beta_3$ for $D_{31}$, $\beta=\beta_3$ for $D_{32}$. All other roots from these rook placements either are maximal or satisfy the conditions of Proposition~\ref{prop:distinct_orbits} for an appropriate simple root $\alpha_0$, as we mentioned above. This completes the proof.
 \end{proof}

\subsection*{Acknowledgements}
This work was performed under the development program of Volga Region Mathematical Center (agreement no. 075--02--2021--1393).


\EditInfo{July 07, 2021}{July 28, 2021}{Ivan Kaygorodov}

\end{paper}